\documentclass[aap,MSNbibl,citesort,dvips]{arximspdf}
\usepackage{mathbh}
\usepackage{dcolumn}
\usepackage{graphicx}

\doi{10.1214/11-AAP790}
\volume{22}
\issue{2}
\pubyear{2012}
\firstpage{860}
\lastpage{880}

\makeatletter

\newcolumntype{d}[1]{D{.}{.}{#1}}

\newtheorem{theor}{Theorem}
\newtheorem{lemma}[theor]{Lemma}

\newcommand{\Z}{\mathbb{Z}}
\newcommand{\ind}{{\mathbh1}}
\newcommand{\ep}{\varepsilon}

\newcommand{\cv}{\operatorname{cv}}
\newcommand{\ax}{\operatorname{ax}}
\newcommand{\vm}{\operatorname{vm}}
\newcommand{\card}{\operatorname{card}}
\newcommand{\bin}{\operatorname{Binomial}}

\makeatother

\begin{document}
\begin{frontmatter}

\title{The Axelrod model for the dissemination of~culture revisited}
\runtitle{\hspace*{-12pt}The Axelrod model for the dissemination of culture revisited}

\begin{aug}
\author[A]{\fnms{Nicolas} \snm{Lanchier}\corref{}\thanksref{t1}\ead[label=e1]{lanchier@math.asu.edu}}
\runauthor{N. Lanchier}
\affiliation{Arizona State University}
\address[A]{School of Mathematical\\
\quad and Statistical Sciences\\
Arizona State University\\
Tempe, Arizona 85287\\
USA} 
\end{aug}

\thankstext{t1}{Supported in part by NSF Grant DMS-10-05282.}

\received{\smonth{4} \syear{2011}}

%
\begin{abstract}
This article is concerned with the Axelrod model, a stochastic process
which similarly to the voter model includes social
influence, but unlike the voter model also accounts for homophily.
Each vertex of the network of interactions is characterized by a set of
$F$ cultural features, each of which can assume
$q$ states.
Pairs of adjacent vertices interact at a rate proportional to the
number of features they share, which results in the
interacting pair having one more cultural feature in common.
The Axelrod model has been extensively studied during the past ten
years, based on numerical simulations and simple mean-field
treatments, while there is a total lack of analytical results for the
spatial model itself.
Simulation results for the one-dimensional system led physicists to
formulate the following conjectures.
When the number of features $F$ and the number of states $q$ both equal
two, or when the number of features exceeds
the number of states, the system converges to a monocultural
equilibrium in the sense that the number of cultural domains
rescaled by the population size converges to zero as the population
goes to infinity.
In contrast, when the number of states exceeds the number of features,
the system freezes in a~highly fragmented configuration
in which the ultimate number of cultural domains scales like the
population size.
In this article, we prove analytically for the one-dimensional system
convergence to a~monocultural equilibrium in terms of
clustering when $F = q = 2$, as well as fixation to a highly fragmented
configuration when the number of states is sufficiently
larger than the number of features.
Our first result also implies clustering of the one-dimensional
constrained voter model.
\end{abstract}

%

\end{frontmatter}

\section{Introduction}
\label{secintro}

Opinion and cultural dynamics are driven by social influence, the
tendency of individuals to become more similar
when they interact,\vadjust{\goodbreak} which is the basic mechanism of the voter model
introduced independently by Clifford and
Sudbury \cite{cliffordsudbury1973} and Holley and Liggett \cite
{holleyliggett1975}.
Social influence alone usually drives the system to a monocultural
equilibrium, whereas differences between individuals and
groups persist in the real world.
In his seminal paper \cite{axelrod1997}, political scientist Robert
Axelrod explains the diversity of cultures
as a consequence of homophily, which is the tendency to interact more
frequently with individuals which are more similar.
In his model, actors are characterized by a finite number of cultural features.
In Axelrod's own words, the more similar an actor is to a neighbor, the
more likely that actor will adopt one of the
neighbor's traits.
The network of interactions is a finite connected graph $G$ with vertex
set $V$ and edge set $E$, where each vertex $x$ is
characterized by a vector $X (x)$ of $F$ cultural features, each of
which assuming $q$ possible states,
\[
X (x) = (X^1 (x), \ldots, X^F (x)),
\]
where
\[
X^i (x) \in\{1, 2, \ldots, q \} \qquad\mbox{for } i = 1, 2, \ldots,
F.
\]
At each time step, a vertex $x$ is picked uniformly at random from the
vertex set along with one of its neighbors $y$.
Then, with a probability equal to the fraction of features $x$ and $y$
have in common, one of the features for which states are
different (if any) is selected, and the state of vertex $x$ is set
equal to the state of vertex $y$ for this cultural feature.
Otherwise nothing happens.
In order to describe more generally the Axelrod dynamics on both finite
and infinite graphs, we assume that the system
evolves in continuous-time with each pair of adjacent vertices
interacting at rate one, which causes one of the two vertices
chosen uniformly at random to mimic the other vertex in the case of an update.
This induces a continuous-time Markov process whose state at time $t$
is a function $X_t$ that maps the vertex set of the
graph into the set of cultures $\{1, 2, \ldots, q \}^F$, and whose
dynamics are described by the generator $\Omega_{\ax}$
defined on the set of cylinder functions by
\[
\Omega_{\ax} f (X) = \sum_{x \in V} \sum_{x \sim y} \sum
_{i = 1}^F \frac{1}{2F}
\biggl[\frac{F (x, y)}{1 - F (x, y)} \biggr] \ind\{X^i (x) \neq X^i
(y) \} [f (X_{y \to x}^i) - f (X)],
\]
where $x \sim y$ means that $x$ and $y$ are connected by an edge,
\[
X_{y \to x}^i (x) = (X^1 (x), \ldots, X^{i - 1} (x), X^i (y), X^{i
+ 1} (x), \ldots, X^F (x))
\]
and $X_{y \to x}^i (z) = X (z)$ for all $z \neq x$, and
\[
F (x, y) = \frac{1}{F} \sum_{i = 1}^F \ind\{X^i (x) = X^i
(y) \}
\]
denotes the fraction of cultural features vertices $x$ and $y$ have in common.
To explain\vadjust{\goodbreak} the expression of the Markov generator of the Axelrod model,
we note that
\[
\frac{1}{2F} \biggl[\frac{F (x, y)}{1 - F (x, y)} \biggr] = F
(x, y) \times\frac{1}{F (1 - F (x, y))} \times\frac{1}{2},
\]
which represents the fraction of features both vertices have in common,
which is the rate at which the vertices interact, times
the reciprocal of the number of features for which both vertices
disagree, which is the probability that any of these features
is the one chosen to be updated, times the probability one half that
vertex $x$ rather than vertex $y$ is chosen to be updated.

The two-feature, two-state Axelrod model is also closely related to the
constrained voter model introduced by
V\'azquez et al. \cite{vazquezkrapivskyredner2003} which identifies
two of the cultures with no common feature
to be a centrist opinion, and the other two cultures to be a leftist
and a rightist opinions.
This results in a~stochastic process somewhat similar to the voter
model except that leftists and rightists are too
incompatible to interact.
Thinking of leftist as $-$ state, centrist as 0 state and rightist as
$+$ state, the constrained voter model is
the Markov process whose state at time $t$ is a function $Z_t$ that
maps the vertex set of the graph into the opinion
set $\{-1, 0, + 1 \}$, and whose dynamics are described by the Markov
generator $\Omega_{\cv}$ defined on the set of cylinder
functions by
\[
\Omega_{\cv} f (Z) = \frac{1}{2} \sum_{x \in V} \sum_{x
\sim y} \sum_{\ep= -1}^1 \ind\{Z (y) = \ep\} \ind\{Z (x) \neq- \ep\} [f
(Z_{x, \ep}) -
f (Z)],
\]
where $Z_{x, \ep}$ is the configuration defined by
\[
Z_{x, \ep} (z) = \ep\qquad\mbox{if } z = x \quad\mbox{and}\quad
Z_{x, \ep} (z) = Z (z) \qquad\mbox{if } z \neq x.
\]
In order to understand these opinion and cultural dynamics, we study
analytically the number and mean size of the cultural
domains at equilibrium.
The number $N (t)$ of cultural domains at time $t$ is the number of
connected components of the graph obtained by removing
all the edges that connect two vertices that do not share the same
culture at time $t$, while the mean size $S (t)$ is defined
as the mean number of vertices per connected component.
Note that when the dynamics take place on a finite connected graph, the
mean size of the cultural domains is also equal to the
total number of vertices divided by the number of cultural domains.

The constrained voter model and Axelrod model have been extensively
studied over the past ten years by social scientists as
well as statistical physicists based on numerical simulations and
simple mean field treatments, while there is a total lack of
analytical results for the spatial models.
We refer the reader to Sections~III.B and IV.A of Castellano et al.
\cite{castellanofortunatoloreto2009} for a review,
and references therein for more details about numerical results.
Because spatial simulations are usually difficult to interpret, there
is a need for rigorous analytical results, and this article is
intended to provide analytical proofs of important conjectures
suggested by the simulations, which also gives insight into
the mechanisms that promote convergence to either a monocultural
equilibrium or, on the contrary, a highly fragmented
configuration where cultural domains are uniformly bounded.\vspace*{8pt}

\textit{Convergence to a monocultural equilibrium.}
Letting $\theta$ denote the initial density of centrists in the
one-dimensional constrained voter model, the mean-field
analysis in \cite{vazquezkrapivskyredner2003} suggests that the
average domain length at equilibrium is
\[
\lim_{t \to\infty} E (S (t)) \sim N^{2 \psi(\theta)}
\qquad\mbox{where }
\psi(\theta) = - \frac{1}{8} + \frac{2}{\pi^2} \biggl[\cos^{-1} \biggl(\frac{1 - 2
\theta}{\sqrt2} \biggr) \biggr]^2,
\]
when the length $N$ of the system is large.
V\'azquez et al. \cite{vazquezkrapivskyredner2003} also showed that
these predictions agree with their numerical simulations
when the initial density of centrists is small enough, as indicated by
their Figure 5, from which they conclude that, for small
$\theta$, the system ends up with high probability to a frozen mixture
of leftists and rightists.
Their simulations, however, also suggest that a typical final frozen
state is characterized by two spatial scales, with few
cultural domains covering macroscopically large fractions of the
universe and a number of small domains.
The presence of two spatial scales also holds for the two-feature
two-state Axelrod model.
Even though it was not the conclusion of V\'azquez et al. \cite
{vazquezkrapivskyredner2003}, this somewhat suggests
convergence to a monocultural equilibrium in which the number of
cultural domains does not scale like the population size, that is,
the cardinality of the vertex set.
Our first result shows that both the constrained voter model and the
two-feature two-state Axelrod model on the one-dimensional
infinite lattice indeed converge to a monocultural equilibrium in terms
of a clustering similar to that of the one-dimensional
voter model.
This clustering indicates that the only stationary distributions are
the ones supported on the set
of configurations in which all vertices share the same culture.
\begin{theor}
\label{thconsensus}
Starting from a translation invariant product measure in which each
culture/opinion occurs with positive probability,
the one-dimensio\-nal two-feature two-state Axelrod model and the
one-dimensional constrained voter model cluster, that is,
\[
\lim_{t \to\infty} P \bigl(X_t (x) \neq X_t (y)\bigr) = \lim_{t \to
\infty} P \bigl(Z_t (x) \neq Z_t (y)\bigr) = 0
\qquad\mbox{for all } x, y \in\Z.
\]
\end{theor}

The apparent contradiction between our analytical result and the
numerical results of \cite{vazquezkrapivskyredner2003} is
due to the fact that, when taking place on a finite connected graph,
the dynamics of the Axelrod model and constrained voter
model may drive the system to a culturally fragmented frozen
configuration even though they promote convergence to a monocultural
equilibrium, which again reveals the difficulty in interpreting spatial
simulations and the need for analytical results.\vspace*{8pt}

\textit{Fixation to a fragmented configuration.}
For the one-dimensional Axelrod model with an arbitrary number of
features and states per feature,
Vilone et al. \cite{vilonevespignanicastellano2002} have predicted
through the analysis of mean-field approximation
supported by simulation results that convergence to a monocultural
equilibrium occurs when $F > q$ whereas fixation to a
highly fragmented configuration occurs when $F < q$.
Our second result establishes partly the latter in the sense that the
expected number of cultural domains at equilibrium on a
path-like graph scales like the length of the graph for an infinite
proper subset of the parameter region $F < q$.
\begin{theor}
\label{thcoexistence}
Assume that $G = \{0, 1, \ldots, N \}$ and $F < q$.
Then, starting from a translation invariant product measure in which
each culture occurs with the same probability $q^{-F}$,
\[
N^{-1} \lim_{t \to\infty} E (N (t)) \geq\biggl(1 - \frac{1}{q} \biggr)^F + \frac
{F}{q - F} \biggl(
\biggl(1 - \frac{1}{q} \biggr)^F - \biggl(1 - \frac{1}{q} \biggr) \biggr).
\]
\end{theor}

Note that the lower bound for the expected number of cultural domains
also gives an upper bound for the expected length of the
cultural domains since the expected number of domains times their
expected length is equal to the cardinal of the vertex set.
Note also that the theorem does not a priori exclude clustering of the
system on the infinite one-dimensional lattice, but
it strongly suggests the latter since it gives upper bounds for the
expected length of the cultural domains which are uniform in
the size of the network of interactions.
Finally, we point out that, even though our estimate
holds for all $F < q$, the
theorem does not fully prove the conjecture of \cite
{vilonevespignanicastellano2002} since the lower bound can be negative:
to make the coexistence region more explicit when both parameters are
large, we fix $F/q = x$ and notice that the lower
bound in Theorem \ref{thcoexistence} converges to
\[
\biggl(\frac{1}{x} - 1 \biggr)^{-1} \biggl(\frac{e^{-x}}{x} - 1 \biggr)
\qquad\mbox{as } F \to\infty.
\]
The equation $e^{-x} = x$ has a unique positive solution given by $x_0
\approx0.567$, which indicates that the one-dimensional
Axelrod model coexists in the sense that the expected number of
cultural domains scales like the population size
when $F < x_0 \times q$ and $F$ is large.
For explicit conditions for coexistence when the parameters are small,
we refer the reader to Table \ref{tabdomains} which gives
some values of the upper bound for the expected length of the cultural
domains at equilibrium.

%
\begin{table}
\tabcolsep=0pt
\caption{Upper bounds for the expected length of the
cultural domains at equilibrium}
\label{tabdomains}
\begin{tabular*}{\tablewidth}{@{\extracolsep{\fill}}lccd{2.4}d{2.4}ccccc@{}}
\hline
& \multicolumn{1}{c}{$\bolds{q = 4}$}
& \multicolumn{1}{c}{$\bolds{q = 8}$}
& \multicolumn{1}{c}{$\bolds{q = 12}$}
& \multicolumn{1}{c}{$\bolds{q = 16}$}
& \multicolumn{1}{c}{$\bolds{q = 20}$}
& \multicolumn{1}{c}{$\bolds{q = 24}$}
& \multicolumn{1}{c}{$\bolds{q = 28}$}
& \multicolumn{1}{c}{$\bolds{q = 32}$}
& \multicolumn{1}{c}{$\bolds{q = 36}$} \\
\hline
$F = 2$ & 2.6667 & 1.3714 & 1.2121 & 1.1487 & 1.1146 & 1.0932 & 1.0785
& 1.0679 & 1.0597 \\ [2pt]
$F = 3$ & neg. & 1.8286 & 1.3861 & 1.2535 & 1.1890 & 1.1508 & 1.1255 &
1.1074 & 1.0940 \\ [2pt]
$F = 4$ & -- & 3.3629 & 1.6645 & 1.3938 & 1.2810 & 1.2188 & 1.1792 &
1.1519 & 1.1318 \\ [2pt]
$F = 5$ & -- & neg. & 2.1989 & 1.5943 & 1.3985 & 1.3007 & 1.2417 &
1.2022 & 1.1738 \\ [2pt]
$F = 6$ & -- & neg. & 3.7048 & 1.9091 & 1.5552 & 1.4017 & 1.3154 &
1.2599 & 1.2211 \\ [2pt]
$F = 7$ & -- & neg. & 45.641 & 2.4851 & 1.7767 & 1.5304 & 1.4040 &
1.3268 & 1.2746 \\ [2pt]
$F = 8$ & -- & -- & \multicolumn{1}{c}{neg.} & 3.9072 & 2.1170 & 1.7007 & 1.5132 & 1.4058
& 1.3360 \\ [2pt]
$F = 9$ & -- & -- & \multicolumn{1}{c}{neg.} & 13.637 & 2.7127 & 1.9385 & 1.6514 & 1.5005
& 1.4071 \\
\hline
\end{tabular*}
\end{table}

The intuition behind Theorems \ref{thconsensus} and \ref
{thcoexistence} that appears in our proofs can be interpreted
in terms of active versus frozen boundaries between adjacent cultural domains.
Here, we call an active boundary the boundary between two cultural
domains with at least one feature in common.
Even though the infinite system keeps evolving indefinitely, Theorem
\ref{thcoexistence} indicates that for reasonably large
values of the number of states $q$ the incompatibility between adjacent
vertices prevents a positive fraction of boundaries to
ever become active, that is, a positive fraction of the boundaries
frozen initially stay frozen at any time.
In contrast, the result of Theorem \ref{thconsensus} is symptomatic of
a~large activity of the system in the sense that each vertex
changes its culture infinitely often which results in the destruction
of the frozen boundaries thus in the presence of cultural
domains that keep growing indefinitely.
Cultural dynamics including two features and two states per feature
taking place on finite graphs operate similarly by promoting
convergence to a monocultural equilibrium.
However, due to the finiteness of the network of interactions, the
system may fixate before reaching a total consensus in which
case the final frozen configuration is characterized by two spatial scales:
cultural domains whose length scales like the size of the system and
domains which are uniformly bounded, as observed
in V\'azquez et al. \cite{vazquezkrapivskyredner2003}.
The rest of the article is devoted to the proofs.


\section{\texorpdfstring{Proof of Theorem \protect\ref{thconsensus}}{Proof of Theorem 1}}
\label{secconsensus}

Note first that the constrained voter model is obtained from the
two-feature two-state Axelrod model by identifying
two cultures without common features with the centrist opinion, and
each of the other two cultures with the leftist and rightist
opinions, respectively.
Therefore the mean cluster size is stochastically larger for the
constrained voter model than the Axelrod model, so it suffices
to prove the result for the latter.
To study the probability of a consensus when $F = q = 2$, but also the
expected number of cultural domains at
equilibrium when $F < q$ in the next section, the idea is to analyze
the evolution of the agreements along the edges rather than
the actual opinion at each vertex.
The network can be viewed as a weighted graph where each edge is
assigned a weight that counts the number of features
its endpoints have in common.
We call $e = \{x, x + 1 \} \in E$ an edge with weight $j$ at time $t$,
or simply a $j$-edge at time $t$, whenever
\[
\bar X_t (e) = \sum_{i = 1}^F \ind\{X_t^i (x) = X_t^i (x + 1)
\} = j.
\]
The key to proving Theorem \ref{thconsensus} is to observe that, when
$F = 2$, clustering of the Axelrod model is equivalent to
almost sure extinction of the 1-edges and the \mbox{0-edges}.
The former follows from the clustering of a certain voter model coupled
with the Axelrod model, while the latter follows from the
combination of the clustering and the site recurrence property of the
voter model that we define below.
Before going into the details of the proof, we start by collecting
important results about the connection between the Axelrod
model, the voter model and coalescing random walks.

%
\begin{figure}

\includegraphics{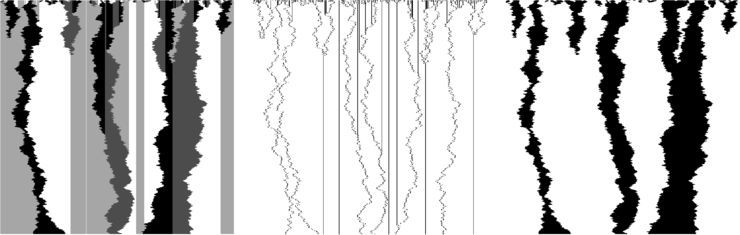}%
\vspace*{-3pt}
\caption{Coupling between the Axelrod model and the voter model.}
\label{figAxM-D1}
\vspace*{-2pt}
\end{figure}

The first ingredient is to observe, as pointed out by V\'azquez and
Redner \cite{vazquezredner2004}, that one recovers the
voter model from the two-feature two-state Axelrod model by identifying
cultures that have no feature in common (see Figure \ref{figAxM-D1} for
simulation pictures of
the coupled Axelrod model and voter model).
Indeed, when $F = q = 2$, we have
\begin{eqnarray*}
\Omega_{\ax} f (X) &=& \frac{1}{4} \sum_{x \in\Z} \sum_{x \sim y} \sum_{i
\neq j}
\ind\{X^i (x) \neq X^i (y) \} \ind\{X^j (x) = X^j (y) \}\\[-4pt]
&&\hphantom{\frac{1}{4} \sum_{x \in\Z} \sum_{x \sim y} \sum_{i
\neq j}}{}
\times [f (X_{y
\to x}^i) - f (X)].
\end{eqnarray*}
Therefore, letting $Y (x) = |X^1 (x) - X^2 (x)|$ for all $x \in\Z$
and noticing that
\[
\{Y (x) \neq Y (y) \} = \bigcup_{i \neq j} \{X^i (x) \neq X^i
(y) \} \cap\{X^j (x) = X^j (y) \},
\]
we obtain that $\{Y_t \dvtx t \geq0 \}$ is the Markov process with generator
\[
\Omega_{\vm} f (Y) = \frac{1}{4} \sum_{x \in\Z} \sum
_{x \sim y} \ind\{Y (x) \neq Y (y) \} [f (Y_{y \to x}) - f
(Y)],\vadjust{\goodbreak}
\]
where $Y_{y \to x}$ is the configuration defined by
\[
Y_{y \to x} (z) = Y (y) \qquad\mbox{if } z = x \quad\mbox{and}\quad
Y_{y \to x} (z) = Y (z) \qquad\mbox{if } z \neq x.
\]
This indicates that $\{Y_t \dvtx t \geq0 \}$ is a time change of the voter
model run at rate $1/2$, but since we are only interested
in the limiting distribution of the Axelrod model, we shall for
simplicity speed up time by a factor two in order to get the usual
voter model run at rate 1.

The voter model can be constructed graphically using an idea of Harris
\cite{harris1972}, which also allows us to exhibit a
duality relationship between the voter model and coalescing random walks.
This construction is now standard, so we only give a brief description.
To each vertex $x \in\Z$, we attach a Poisson process with parameter one.
Then, at the arrival times of this process, we choose one of the two
neighbors $x \pm1$ uniformly at random, then draw an arrow
from vertex $x \pm1$ to vertex $x$ and put a $\delta$ at $x$ to
indicate that $x$ updates its opinion by mimicking $x \pm1$.
The connection between the voter model and coalescing random walks
appears when keeping track of the ancestry of each vertex going
backward in time which also defines the so-called dual process.
We say that there is a dual path from $(x, T)$ to $(y, T - s)$ if there
are sequences of times and vertices
\[
s_0 = T - s < s_1 < \cdots< s_{n + 1} = T
\quad\mbox{and}\quad
x_0 = y, x_1, \ldots, x_n = x
\]
such that the following two conditions hold:
\begin{longlist}[(2)]
\item[(1)] for $i = 1, 2, \ldots, n$, there is an arrow from $x_{i -
1}$ to $x_i$ at time $s_i$;
\item[(2)] for $i = 0, 1, \ldots, n$, the vertical segment $\{x_i \}
\times(s_i, s_{i + 1})$ does not contain any $\delta$'s.
\end{longlist}
Then, for $A \subset\Z$ finite, the dual process starting at $(A, T)$
is the set-valued process
\begin{eqnarray*}
\hat Y_s (A, T) &=& \{y \in\Z\mbox{: there is a dual path}\\
&&\hphantom{\{}\hspace*{0pt}\mbox{from
$(x, T)$ to $(y, T - s)$ for some $x \in A$} \}.
\end{eqnarray*}
The dual process is naturally defined only for dual times $0 \leq s
\leq T$.
However, it is convenient to assume that the Poisson processes in the
graphical representation are also defined for
negative times so that the dual process can be defined for all $s \geq0$.
The reason for introducing the dual process is that it allows one to
deduce the state of the process at the current time from the
configuration at earlier times based on the duality relationship
\[
Y_t (x) = Y_{t - s} (\hat Y_s (x, t)) = Y_0 (\hat Y_t (x, t))
\qquad\mbox{for all } s \in(0, t).
\]
Moreover, it can be seen from the graphical representation that the
dual process evolves according to a system of simple symmetric
coalescing random walks that jump at rate 1, so questions about the
voter model can be answered by looking at this system of
coalescing random walks.

The first step, as previously mentioned, is to establish extinction of
the 1-edges, which is directly related to the
clustering of the one-dimensional voter model.
\begin{lemma}
\label{lem1-edges}
There is almost sure extinction of the 1-edges, that is,
\[
\lim_{t \to\infty} P \bigl(\bar X_t (e) = 1\bigr) = 0 \qquad\mbox{for
all } e \in E.
\]
\end{lemma}

\begin{pf}
Since the one-dimensional voter model clusters \cite
{cliffordsudbury1973,holleyliggett1975}, we have
\begin{eqnarray*}
&&\lim_{t \to\infty} P \bigl(\bar X_t (\{x, x + 1 \}) =
1\bigr)\\
&&\qquad =
\lim_{t \to\infty} P \bigl(X_t^1 (x) \neq X_t^1 (x + 1)
\mbox{ and } X_t^2 (x) = X_t^2 (x + 1)\bigr) \\
&&\qquad\quad{} +
\lim_{t \to\infty} P \bigl(X_t^1 (x) = X_t^1 (x + 1) \mbox{ and } X_t^2 (x)
\neq X_t^2 (x + 1)\bigr) \\
&&\qquad =
\lim_{t \to\infty} P \bigl(Y_t (x) \neq Y_t (x + 1)\bigr) = 0
\end{eqnarray*}
for every vertex $x \in\Z$, which proves extinction of the 1-edges.
In words, the \mbox{1-edges} in the Axelrod model correspond to the interfaces
of the underlying voter model, which evolve according to
a system of annihilating random walks that, because of clustering, goes extinct.
\end{pf}

The second step is to prove that there is almost sure extinction of the
0-edges, which follows from the combination of clustering and
the following property:
\[
P \bigl(Y_t (x) \neq Y_s (x) \mbox{ for some } t > s\bigr) = 1
\qquad\mbox{for all } x \in\Z\mbox{ and } s > 0
\]
that we shall call site recurrence of the one-dimensional voter model.
The terminology is motivated by the article of Erd\H os and Ney \cite
{erdosney1974} who conjectured that, given a system of discrete-time
annihilating random walks starting with one particle at each site
except the origin, the probability that the origin is visited infinitely
often is one, property that Arratia \cite{arratia1983} called later
site recurrence.
The continuous-time version of their conjecture has been proved by
Schwartz \cite{schwartz1978} based on the connection with the
one-dimensional voter model whose interfaces precisely evolve according
to a system of annihilating random walks.
In particular, the idea of her proof is to show that the
one-dimensional voter model starting from a particular deterministic
configuration
is site recurrent in the sense defined above.
Her result easily extends to the process starting from more general
random configurations, but we give a somewhat shorter proof in
Lemma \ref{lemrecurrence} below.
Note that the site recurrence of the voter model in higher dimensions
directly follows from the law of large numbers for the occupation time of
the process established in Cox and Griffeath \cite{coxgriffeath1983}
which does not hold in one dimension.
\begin{lemma}
\label{lemrecurrence}
The one-dimensional voter model is site recurrent, that is,
\[
P \bigl(Y_t (x) \neq Y_s (x) \mbox{ for some } t > s\bigr) = 1
\qquad\mbox{for all } x \in\Z\mbox{ and } s > 0.
\]
\end{lemma}

\begin{pf}
The key is to observe that, for all $y \in\Z$, the process that keeps
track of the number of vertices at time $t$ that descend from $y$; namely
\[
M_t (y) := \card\{z \in\Z\mbox{: there is a dual path from
$(z, t)$ to $(y, 0)$} \}
\]
is a martingale absorbed at state 0.
Note that, though stated for the one-dimensional voter model only, this
property holds in any spatial dimension.
Since in addition $M_t (y)$ is an integer-valued process, a
straightforward application of the martingale convergence theorem implies
that it converges almost surely to its absorbing state.
Now, let $x \in\Z$ and $s > 0$, and define
\[
\Phi(s) := \inf\{t > 0 \dvtx M_t (y) = 0 \} \qquad\mbox{where }
y := \hat Y_s (x, s).
\]
Since there is a dual path from $(x, s)$ to $(y, 0)$, this stopping
time is larger than~$s$, but almost sure convergence of $M_t (y)$ to
zero implies that time $\Phi(s)$ is almost surely finite.
In addition,
\[
\mbox{there is no dual path from $(x, \Phi(s))$ to $(y, 0)$}
\]
from which it follows that the spin at $(x, \Phi(s))$ and the spin at
$(x, s)$ originate from different vertices at time 0 and thus
are independent since the initial configuration is distributed
according to a product measure.
This holds for all $s > 0$ hence defining recursively $s_0 = s$ and
\[
s_{i + 1} := \Phi(s_i) = \inf\{t > 0 \dvtx M_t (y) = 0 \}
\qquad\mbox{where } y := \hat Y_{s_i} (x, s_i)
\]
for all integers $i \geq0$ induces an increasing sequence of stopping
times which are all almost surely finite.
Moreover, the collection of spins at $(x, s_i)$ are independent,
determined from the spins of different vertices at time 0.
It is straightforward to deduce that
\[
P \bigl(Y_s (x) := Y_{s_0} (x) = Y_{s_1} (x) = \cdots= Y_{s_i} (x)\bigr) \to
0 \qquad\mbox{as } i \to\infty
\]
since each type occurs initially with positive probability.
The lemma follows.
\end{pf}

\begin{lemma}
\label{lem0-edges}
There is almost sure extinction of the 0-edges, that is,
\[
\lim_{t \to\infty} P \bigl(\bar X_t (e) = 0\bigr) = 0 \qquad\mbox{for
all } e \in E.
\]
\end{lemma}

\begin{pf}
First of all, we observe that, since the initial configuration, as well
as the evolution rules of the process, are translation invariant
in space, the probability of edge $e$ being a 0-edge at time $t$ does
not depend on the specific choice of $e$.
This key property is implicitly used repeatedly in the proof of the lemma.
Now, we let $0 < s < t < \infty$, and partition the set of 0-edges at
time $t$ into the subset $\Omega_-$ of those edges that have been lately
updated by time $s$ and the subset $\Omega_+$ of those edges that have
been lately updated after time $s$, namely
\begin{eqnarray*}
\Omega_- & = & \{e \in E \dvtx\bar X_u (e) = 0 \mbox{ for all } u \in
(s, t) \}, \\
\Omega_+ & = & \{e \in E \dvtx\bar X_t (e) = 0 \mbox{ and } \bar X_u
(e) = 1 \mbox{ for some } u \in(s, t) \}.
\end{eqnarray*}
Note that an update at vertex $x$ in the Axelrod model corresponds to
the simultaneous update of the pair of edges incident to this vertex.
In addition, the culture at $x$ flips at a positive rate if and only if
$x$ has exactly one feature in common with one of its
two nearest neighbors, that is, if and only if at least one of both
edges incident to $x$ is a 1-edge.
Also, since only one feature changes at a time, when an edge pair is
updated, the weight of each edges varies by exactly one unit.
In particular, accounting for symmetry, there are only four possible
transitions of the edge pairs:
\[
(1, 0) \to(2, 1),\qquad (1, 1) \to(2, 0),\qquad (1, 1) \to(2, 2),\qquad
(1, 2) \to(2, 1).
\]
It follows that the probability of $e$ being a 1-edge is nonincreasing
(this can be seen from the fact that the 1-edges evolve according to
a system of annihilating random walks) and that a 0-edge can only
result from the annihilation of two 1-edges.
This, together with Lemma \ref{lem1-edges} which claims extinction of
the 1-edges, implies that, for all $\ep> 0$, there exists $s$ large
such that
\[
2 \times P (e \in\Omega_+) \leq P \bigl(\bar X_s (e) = 1\bigr) \leq\ep.
\]
Time $s$ being fixed, Lemma \ref{lemrecurrence} implies the existence
of $t > s$ such that
\begin{eqnarray*}
P (e = \{x, x + 1 \} \in\Omega_-) & = &
P \bigl(X_u^i (x) \neq X_u^i (x + 1) \mbox{ for all } u
\in(s, t) \mbox{ and } i = 1, 2\bigr) \\
& \leq&
P \bigl(Y_u (x) = Y_u (x + 1) \mbox{ for all } u \in(s,
t)\bigr) \\
& \leq&
P \bigl(Y_u (x) = Y_s (x) \mbox{ for all } u \in(s, t)\bigr)
\leq\ep.
\end{eqnarray*}
Combining the previous two estimates, we obtain that, for all $\ep> 0$
small, there exists a large but finite time $s > 0$ and a
large but finite time $t > s$ such that
\[
P \bigl(\bar X_t (e) = 0\bigr) = P (e \in\Omega_+) + P (e \in
\Omega_-) \leq2 \ep,
\]
which establishes extinction of the 0-edges.
\end{pf}

Having established that both sets of 0-edges and 1-edges go extinct,
the proof of Theorem \ref{thconsensus} is now straightforward,
and follows the lines of Lemma \ref{lem1-edges}.
While the latter shows that clustering of the voter model implies
extinction of the 1-edges, the last step is to prove that,
conversely, extinction of type 0 and type 1 edges implies clustering of
the Axelrod model.
Fix $x < y$, and let
\[
z_0 = x < z_1 < \cdots< z_k = y \qquad\mbox{with }
k = |x - y|.
\]
Denote by $e_i = \{z_i, z_{i + 1} \}$ the edge connecting vertex $z_i$
and vertex $z_{i + 1}$.
Then, extinction of both the 1-edges and the 0-edges given,
respectively,
by Lemmas \ref{lem1-edges} and \ref{lem0-edges} implies that
\begin{eqnarray*}
&&\lim_{t \to\infty} P \bigl(X_t (x) \neq X_t (y)\bigr) \\
&&\qquad \leq
\lim_{t \to\infty} P \bigl(X_t (z_i) \neq X_t (z_{i +
1}) \mbox{ for some } i = 0, 1, \ldots, k - 1\bigr) \\
&&\qquad\leq
\lim_{t \to\infty} \sum_{i = 0}^{k - 1} P \bigl(\bar
X_t (e_i) = 0 \mbox{ or } \bar X_t (e_i) = 1\bigr) = 0.
\end{eqnarray*}
This completes the proof of Theorem \ref{thconsensus}.


\section{\texorpdfstring{Proof of Theorem \protect\ref{thcoexistence}}{Proof of Theorem 2}}
\label{seccoexistence}

This section is devoted to the proof of Theorem~\ref{thcoexistence}
which again relies on the analysis of the agreements
along the edges rather than the actual opinion at each vertex.
Before going into the details of the proof, we briefly introduce its
main steps.
First, it is noted in Lemma \ref{lemdomains} that the process on a
finite graph reaches almost surely one of its absorbing states
which consist of the configurations in which each edge has either
weight zero or weight $F$.
The ultimate number of cultural domains on a path-like graph is roughly
equal to the ultimate number of edges with weight zero,
so the strategy is to bound from below the number of such edges.
To do so, we introduce
\[
W (t) := \sum_{j = 0}^F j W_j (t) \qquad\mbox{where }
W_j (t) := \card\{e \in E \dvtx\bar X_t (e) = j \}
\]
and where $\bar X_t (e)$ is defined as in the previous section.
In words, $W (t)$ keeps track of the total number of agreements in the system.
Let $y$ be a vertex with degree 2, and let $x$ and $z$ be its two
nearest neighbors.
Then, we observe that if the individual at vertex $y$ updates its
culture at time $t$ by mimicking the $i$th feature of vertex $x$,
then we have the following alternatives:
\begin{longlist}[(2)]
\item[(0)] $W (t) - W (t-) = 0$ whenever $X_{t-}^i (y) = X_{t-}^i
(z)$.
\item[(1)] $W (t) - W (t-) = 1$ whenever $X_{t-}^i (x) \neq X_{t-}^i
(z)$ and $X_{t-}^i (y) \neq X_{t-}^i (z)$.
\item[(2)] $W (t) - W (t-) = 2$ whenever $X_{t-}^i (x) = X_{t-}^i (z)$.
\end{longlist}
The first key step is to prove that if the $i$th feature of vertex $y$
differs from the $i$th feature of both of its
neighbors (cases 1 and 2 above), then the $i$th feature of vertex $x$
and the $i$th feature of vertex $z$ are independent, which is
established in Lemma~\ref{lemindependent}.
It follows that
\[
P \bigl(W (t) - W (t-) = 2 | W (t) \neq W (t-)\bigr) \leq(q - 1)^{-1}
\]
as stated in Lemma \ref{lemmoving}.
Note that two nearest neighbors, say $y$ and $z$, cannot interact as
long as they are connected by an edge with weight zero;
therefore edges with weight zero can possibly change their state only
when two agreements emerge simultaneously (case 2 above).
This, together with the previous inequality, leads us to consider the
following urn problem which is also based on an idea
initially introduced in \cite{lanchier2010}.
There are $F + 1$ boxes labeled from box 0 to box $F$ containing all
together a total of $N$ balls, which corresponds to the number
of edges.
The game starts with as many balls in box $j$ as there are edges with
weight $j$ in the Axelrod model at time~0, and evolves in
discrete time according to the following stochastic rules:
\begin{longlist}[(2)]
\item[(1)] At each time step, we move a ball from box $j$ to box $j +
1$ where (if it exists) box $j$ is chosen uniformly at random from
the set of nonempty inner boxes.
\item[(2)] In case a ball has indeed been moved, and box 0 is
nonempty,\vspace*{1pt} we move an additional ball from box 0 to box 1 only with
probability $(q - 1)^{-1}$.
\item[(3)] The game halts when all the inner boxes are empty.
\end{longlist}
Here, inner boxes refer to boxes $1, 2, \ldots, F - 1$.
Using coupling arguments and the inequality established in Lemma \ref
{lemmoving}, we prove that the expected number of balls in box 0
when the game halts is larger than the expected number of edges with
weight zero in the Axelrod model when the system gets trapped.
This is done in Lemma \ref{lemballs}.
Finally, the lower bound of Theorem \ref{thcoexistence} is proved to
be a lower bound for the expected number of balls in box 0
when the game halts, and thus a lower bound for the expected value of
the ultimate number of domains, in Lemmas \ref{leminitial} and
\ref{lemfinal}.
The key to keeping track of the number of balls in box 0 is to divide
the game into rounds and paint the balls with two different
colors at the beginning of each round.
\begin{lemma}
\label{lemdomains}
We have $\lim_{ t \to\infty} N (t) - W_0 (t) = 1$.
\end{lemma}
\begin{pf}
The number $N (t)$ of cultural domains at time $t$ is the number of
connected components of the graph obtained by removing all
edges whose weight at time $t$ differs from $F$.
In the case of a finite tree, this results in a forest whose number of
connected components is equal to the number of edges
removed plus one, from which it follows that
\[
N (t) = W_0 (t) + W_1 (t) + \cdots+ W_{F - 1} (t) +
1.
\]
Since the Axelrod model on a finite graph converges to one of its
absorbing states and that each absorbing state is characterized
by all the edges having weight either zero or $F$, we also have
\[
\lim_{t \to\infty} W_j (t) = 0 \qquad\mbox{for } j = 1, 2,
\ldots, F - 1.
\]
The result follows.
\end{pf}
\begin{lemma}
\label{lemindependent}
Let $0 \leq x < y < z \leq N$ and fix $i \in\{1, 2, \ldots, F \}$. Then
\[
P \bigl(X_t^i (x) = X_t^i (z) | X_t^i (x) \neq X_t^i (y) \mbox{ and }
X_t^i (y) \neq X_t^i (z)\bigr) = (q - 1)^{-1}.
\]
\end{lemma}

\begin{pf}
The idea is that, given that the pairs $\{x, y \}$ and $\{y, z \}$
disagree on their $i$th feature, the states of this feature for $x$ and $z$
must originate from different vertices at time 0.
To make this argument rigorous, we first construct the Axelrod model
graphically from collections of independent random variables:
for each oriented edge $e = (u, v)$ and $n \geq1$:
\begin{longlist}[(2)]
\item[(1)] $T_n (e)$ is the $n$th arrival time of a Poisson process
with rate $1/2$;
\item[(2)] $U_n (e)$ is the discrete random variable uniformly
distributed over $\{1,\allowbreak 2, \ldots, F \}$;
\item[(3)] $V_n (e) = \{V_{n, k} (e) \dvtx k \geq1 \}$ is an infinite
sequence of independent random variables which are uniformly
distributed over the
set of features $\{1, 2, \ldots, F \}$.
\end{longlist}
The process starting from any initial configuration is constructed
inductively as follows.
Assume that the process has been constructed up to time $t-$ where $t =
T_n (e)$, and let
\[
I_{t-} (e) = I_{t-} ((u, v)) = \{i \dvtx X_{t-}^i (u) = X_{t-}^i (v) \}
\]
and
\[
J_{t-} (e) = \{i \dvtx X_{t-}^i (u) \neq X_{t-}^i (v) \}.
\]
In words, the sets $I_{t-} (e)$ and $J_{t-} (e)$ are the sets of
features for which the vertices $u$ and $v$ agree and disagree, respectively.
Then we have the following alternative:
\begin{longlist}[(2)]
\item[(1)] if $U_n (e) \in I_{t-} (e)$ and $J_{t-} (e) \neq
\varnothing$, then we draw an arrow from $u$ to $v$ at time $t$;
\item[(2)] if $U_n (e) \in J_{t-} (e)$ or $J_{t-} (e) = \varnothing
$, then we do nothing.
\end{longlist}
Thinking of an arrow oriented from $u$ to $v$ as representing an
interaction that causes vertex $v$ to mimic one of the features
of $u$, and noticing that
\[
P \bigl(U_n (e) \in I_{t-} (e)\bigr) = \frac{1}{F} \card(I_{t-} (e))
= \frac{1}{F} \sum_{i = 1}^F \ind\{X_{t-}^i (u) = X_{t-}^i
(v) \},
\]
this indicates that, as required, adjacent vertices with different
cultures interact at a rate proportional to the number of features they share.
To complete the construction, the last step is to use the random
variables $V_{n, k} (e)$ to determine the feature vertex $v$
mimics on the event that 1 above occurs.
To choose this feature uniformly at random from $J_{t-} (e)$, as
required, we let
\[
m = \inf\{k \geq1 \dvtx V_{n, k} (e) \in J_{t-} (e) \} \quad\mbox
{and}\quad i = V_{n, m} (e)
\]
and label the arrow with an $i$ to indicate that the $i$th feature of
$v$ at time $t$ is set equal to the $i$th feature of vertex $u$.
Since the network of interactions is finite, the times of the Poisson
events can be ordered; therefore the Axelrod model can be
constructed going forward in time using the previous rules.
However, an\vadjust{\goodbreak} idea of Harris \cite{harris1972} allows us as well to
construct the process on infinite lattices in the same manner.
Given an initial configuration in which features are independent and
uniformly distributed, and the collections of independent random variables
introduced above, we draw the arrows along with their label up to time
$t$ following the rules previously described.
Then we say that there is an $i$-lineage from $(u, t)$ to $(w, t - s)$
if there are
\[
s_0 = t - s < s_1 < \cdots< s_{n + 1} = t
\quad\mbox{and}\quad
u_0 = w, u_1, \ldots, u_n = u
\]
such that the following two conditions hold:
\begin{longlist}[(2)]
\item[(1)] for $j = 1, 2, \ldots, n$, there is an $i$-arrow from
$u_{j - 1}$ to $u_j$ at time $s_j$;
\item[(2)] for $j = 0, 1, \ldots, n$, the segment $\{u_j \} \times
(s_j, s_{j + 1})$ does not contain any tips of $i$-arrow.
\end{longlist}
Note that for all $s \in(0, t)$, there exists a unique vertex $w$ such
that (1) and (2) hold.
We define the process that keeps track of the unique $i$-lineage
starting at $(u, t)$ by letting
\[
\hat X_s^i (u, t) = \{w \mbox{: there is an $i$-lineage from $(u,
t)$ to $(w, t - s)$} \},
\]
and refer the reader to Figure \ref{figdomains} for an illustration.
Although $i$-lineages in the Axelrod model are somewhat reminiscent of
%
%
\begin{figure}

\includegraphics{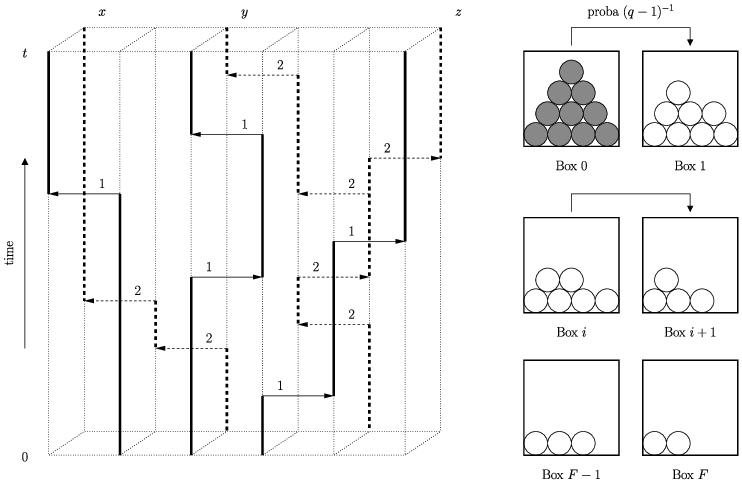}

\caption{Pictures related to the proof of Theorem
\protect\ref{thcoexistence}---lineages in the Axelrod model on the
left-hand side,
and schematic illustration of a single step evolution of the urn
problem on the right-hand side.}
\label{figdomains}
\end{figure}
dual paths in the voter model, the latter can be constructed from
the graphical representation regardless of the initial configuration,
whereas the construction of the former also depends on the initial
configuration, since one needs to construct the process forward in time
up to time $t$ in order to identify the labels on the arrows.
In particular, the culture of a given space--time point cannot be
determined from the initial configuration by simply going backward
in time along the graphical representation.
However, by construction, one has
\[
X_t^i (u) = X_{t - s}^i (\hat X_s^i (u, t)) = X_0^i (\hat X_t^i
(u, t)) \qquad\mbox{for all } 0 \leq s \leq t
\]
from which it follows that
\[
X_t^i (y) \notin\{X_t^i (x), X_t^i (z) \} \quad\mbox{implies}\quad
\hat X_t^i (y, t) \notin\{\hat X_t^i (x, t), \hat X_t^i (z, t) \}.
\]
Moreover, after the random labels on the arrows have been determined by
constructing the process up to time $t$, the system of $i$-lineages
is constructed from the set of $i$-arrows in the same manner as the
system of dual paths in the voter model.
Due, in addition, to the presence of one-dimensional nearest neighbor
interactions, it is straightforward to deduce that, regardless of the
random configuration of the $i$-arrows, the dynamics preserve the order
of the $i$-lineages in the sense that
\[
\hat X_s^i (x, t) \leq\hat X_s^i (y, t) \leq\hat X_s^i (z, t)
\qquad\mbox{for all } 0 \leq s \leq t.
\]
Combining the previous two properties, we deduce that
\[
X_t^i (y) \notin\{X_t^i (x), X_t^i (z) \} \quad\mbox{implies}\quad
\hat X_t^i (x, t) < \hat X_t^i (y, t) < \hat X_t^i (z, t),
\]
indicating that the $i$th feature of vertex $x$ and the $i$th feature
of vertex $z$ at time $t$ are determined by the initial $i$th features of
two different vertices.
Since these features are independent and uniformly distributed at time
0, we obtain that, for all $q_1, q_2 \neq X_t^i (y)$,
\begin{eqnarray*}
&& P \bigl(X_t^i (x) = q_1 \mbox{ and } X_t^i (z) = q_2 | X_t^i (x)
\neq X_t^i (y) \mbox{ and } X_t^i (y) \neq X_t^i (z)\bigr)
\\
&&\qquad=
P \bigl(X_0^i (\hat X_t^i (x, t)) = q_1 \mbox{ and } X_0^i (\hat X_t^i
(z, t)) = q_2 | \hat X_t^i (x, t) \neq\hat X_t^i (z, t) \\
&&\hspace*{142.7pt}\mbox{and } X_t^i (x) \neq
X_t^i (y) \mbox{ and } X_t^i (y) \neq X_t^i (z)\bigr) \\
&&\qquad= (q - 1)^{-2}.
\end{eqnarray*}
The lemma follows by summing over all the $q - 1$ possible values of
$q_1 = q_2$.
\end{pf}

\begin{lemma}
\label{lemmoving}
We have $P (W (t) - W (t-) = 2 | W (t) \neq W (t-)) \leq(q - 1)^{-1}$.
\end{lemma}

\begin{pf}
First, we observe that there exist $n \geq1$ and $e = (x, y)$ such
that $t = T_n (e)$ on the conditional event that the total number of
agreements increases at time $t$.
Then we construct the process up to time $t-$ using the collection of
independent random variables and following the rules introduced
in Lemma \ref{lemindependent} in order to identify the set
\[
J_{t-} (e) = \{1 \leq i \leq F \dvtx X_{t-}^i (x) \neq X_{t-}^i (y)
\},
\]
which is nonempty on the event that $W (t) \neq W (t-)$.
Note that if vertex $y$ is at the boundary of the system, then only the
weight of edge $e$ is updated at time $t$; therefore
\[
P \bigl(W (t) - W (t-) = 2 | W (t) \neq W (t-)\bigr) = 0 \qquad\mbox
{whenever } y = 0 \mbox{ or } y = N.
\]
To deal with the nontrivial case when vertex $y$ has degree 2, observe that
\begin{eqnarray*}
\{W (t) - W (t-) = 2 \} \cap\{V_{n, m} (e) = i \} & = & \{X_{t-}^i
(x) = X_{t-}^i (z) \} \cap\{V_{n, m} (e) = i \}, \\
\{W (t) \neq W (t-) \} \cap\{V_{n, m} (e) = i \} & = & \{X_{t-}^i
(y) \neq X_{t-}^i (z) \} \cap\{V_{n, m} (e) = i \},
\end{eqnarray*}
where vertex $z$ is the unique nearest neighbor of $y$ different from
$x$ and where~$m$ is defined in the construction of the process
given in Lemma \ref{lemindependent}.
Since the random variables $V_{n, k}$, $k \geq1$, are independent of
the configuration of the system at time~$t-$, we have in addition that
\begin{eqnarray*}
&&P \bigl(X_{t-}^i (x) = X_{t-}^i (z) | X_{t-}^i (y) \neq X_{t-}^i (z)
\mbox{ and } V_{n, m} (e) = i\bigr)
\\
&&\qquad=
P \bigl(X_{t-}^i (x) = X_{t-}^i (z) | X_{t-}^i (y) \neq X_{t-}^i (z) \mbox
{ and } i \in J_{t-} (e)\bigr).
\end{eqnarray*}
Combining the previous two properties, we deduce that
\begin{eqnarray*}
&&P \bigl(W (t) - W (t-) = 2 | W (t) \neq W (t-) \mbox{ and } V_{n,
m} (e) = i\bigr) \\
&&\qquad=
P \bigl(X_{t-}^i (x) = X_{t-}^i (z) | W (t) \neq W (t-) \mbox{ and }
V_{n, m} (e) = i\bigr) \\
&&\qquad=
P \bigl(X_{t-}^i (x) = X_{t-}^i (z) | X_{t-}^i (y) \neq X_{t-}^i (z) \mbox
{ and } V_{n, m} (e) = i\bigr)
\\
&&\qquad= P \bigl(X_{t-}^i (x) = X_{t-}^i (z) | X_{t-}^i (y) \neq X_{t-}^i (z) \mbox
{ and } i \in J_{t-} (e)\bigr) \\
&&\qquad= P \bigl(X_{t-}^i (x) = X_{t-}^i (z) | X_{t-}^i (x) \neq X_{t-}^i (y) \mbox
{ and } X_{t-}^i (y) \neq X_{t-}^i (z)\bigr) \\
&&\qquad= (q - 1)^{-1},
\end{eqnarray*}
where the last equality follows from Lemma \ref{lemindependent}.
This completes the proof.
\end{pf}

As previously mentioned, to find a lower bound for the expected value
of the ultimate number of edges with weight zero, we need to
compare the number of such edges with the ultimate number of balls in
box 0 for the game described at the beginning of this section.
To do so, the first step is to couple the Axelrod dynamics with another
urn problem that evolves in continuous-time.
As previously, we start with as many balls in box $j$ as there are
edges with weight $j$ in the Axelrod model at time 0, but the
evolution is now coupled with the cultural dynamics as follows:
\begin{longlist}[(2)]
\item[(0)] if $W (t) - W (t-) = 0$, we do nothing;\vadjust{\goodbreak}
\item[(1)] if $W (t) - W (t-) = 1$, we move a ball from box $j$ to box
$j + 1$ where (if it exists), box $j$ is chosen uniformly at random
from the
set of nonempty inner boxes;
\item[(2)] if $W (t) - W (t-) = 2$, we repeat the same as in 1 above,
and, in case a ball has indeed been moved and box 0 is nonempty, we move
another ball from box 0 to box 1.
\end{longlist}
The game halts when the Axelrod model hits an absorbing state, that is,
when all the edges have either weight zero or weight $F$.
Let $\mathfrak B_j (t)$ denote the number of balls in box $j$ at time $t$.
The next lemma indicates that, at any time, the number of balls in box
0 is smaller than the number of edges with weight zero in the
Axelrod model.
\begin{lemma}
\label{lemballs}
For all $t \geq0$, we have $\mathfrak B_0 (t) \leq W_0 (t)$.
\end{lemma}
\begin{pf}
The intuition behind the result is that a ball is removed from box 0 if
and only if two agreements emerge simultaneously in the
Axelrod model, whereas the latter is only a necessary condition for an
edge with weight zero to change its weight.
To make this argument rigorous, we introduce the key variable $\bar B
(t)$ that represents the number of steps required after time $t$
to move all the balls to box~$F$, excluding the ones which are in box 0
at time $t$, along with its analog for the Axelrod model that
we denote by $\bar W (t)$.
More precisely, we introduce
\[
\bar B (t) := \sum_{j = 1}^F (F - j) \mathfrak B_j (t)
\quad\mbox{and}\quad
\bar W (t) := \sum_{j = 1}^F (F - j) W_j (t).
\]
Then, the idea is to prove by induction that, as long as box 0 is
nonempty (note that once it is empty the result is trivial), we have
\[
\mathfrak B_0 (t) \leq W_0 (t) \quad\mbox{and}\quad \bar B (t) \geq\bar W (t).
\]
The two inequalities to be proved are obviously true at time 0 since
initially there are as many balls in box $j$ as there are
edges with weight $j$.
Assume that they are true at time $t-$ and that a culture is updated at
time $t$. Since
\[
\bar W (t) = F \times\bigl(N - W_0 (t)\bigr) - \sum_{j = 1}^F j W_j
(t) = F \times\bigl(N - W_0 (t)\bigr) - W (t),
\]
and a weight jumps from 0 to 1 at time $t$ only if $W (t) - W (t-) =
2$, we~have:
\begin{longlist}[(2)]
\item[(1)] Assume that $W (t) - W (t-) = 0$. Then
\begin{eqnarray*}
\mathfrak B_0 (t) &=& \mathfrak B_0 (t-) \leq W_0 (t-) = W_0
(t), \\
\bar B (t) &=& \bar B (t-) \geq \bar W (t-) = F \times\bigl(N -
W_0 (t)\bigr) - W (t) = \bar W (t).
\end{eqnarray*}
\item[(2)] Assume that $W (t) - W (t-) = 1$.
Then $\bar W (t-) > 0$ and so $\bar B (t-) > 0$ by assumption.
In particular,\vadjust{\goodbreak} one of the inner boxes is nonempty, which implies that a
ball is indeed moved from some box $j$ to box $j + 1$.
It follows that
\begin{eqnarray*}
\mathfrak B_0 (t) &=& \mathfrak B_0 (t-) \leq W_0 (t-) = W_0
(t), \\
\bar B (t) &=& \bar B (t-) - 1 \geq \bar W (t-) - 1 = F
\times\bigl(N - W_0 (t-)\bigr) - W (t-) - 1 \\
&=& F \times\bigl(N - W_0 (t)\bigr) - W (t) = \bar W (t).
\end{eqnarray*}
\item[(3)] Assume that\vspace*{1pt} $W (t) - W (t-) = 2$.
In case box 0 is empty at time $t-$, the result is trivial.
Otherwise, using as previously that $\bar B (t-) \geq\bar W (t-) > 0$,
we obtain
\begin{eqnarray*}
\mathfrak B_0 (t) &=& \mathfrak B_0 (t-) - 1 \leq W_0 (t-) - 1 \leq
W_0 (t), \\
\bar B (t) &=& \bar B (t-) - 1 + (F - 1) \geq \bar W (t-) - 1 +
(F - 1) \\
& = & F \times\bigl(N - W_0 (t-)\bigr) - W (t) + F \\
& \geq& F \times\bigl(N - W_0 (t) - 1\bigr) - W (t) + F = \bar W (t).
\end{eqnarray*}
\end{longlist}
This completes the proof.
\end{pf}
\begin{lemma}
\label{leminitial}
For all $j = 0, 1, \ldots, F$, we have
\[
E (W_j (0)) = N p_j := N \pmatrix{F \cr j} q^{-j} (1 -
q^{-1})^{F - j}.
\]
\end{lemma}

\begin{pf}
Since initially nearest neighbors agree on each of their features with
probability $q^{-1}$ and that all features are independent,
the probability that a given edge is a $j$-edge is
\[
P (X = j) \qquad\mbox{where } X \sim\bin(F, q^{-1}).
\]
Since the graph contains $N$ edges, the result follows.
\end{pf}
\begin{lemma}
\label{lemfinal}
Assume that $F < q$. Then
\[
N^{-1} \lim_{t \to\infty} E (W_0 (t)) \geq\biggl(1 - \frac{1}{q} \biggr)^F + \frac
{F}{q - F} \biggl(
\biggl(1 - \frac{1}{q} \biggr)^F - \biggl(1 - \frac{1}{q} \biggr) \biggr).
\]
\end{lemma}

\begin{pf}
In view of Lemma \ref{lemballs}, it suffices to prove the result for
$\mathfrak B_0 (t)$ instead of $W_0 (t)$.
First, we consider the discrete-time game introduced at the beginning
of this section.
To count the number of balls more easily, we divide the evolution of
the latter into rounds as follows.
\begin{longlist}[Round 1.]
\item[Round 1.] We paint in black all the balls in box 0 and in
white all the other balls and, at each time step, move a white ball from
box $j$ to box $j + 1$ where (if it exists) box $j$ is chosen uniformly
at random from the set of inner boxes that contain at least one white ball.
In case a white ball has indeed been moved, and box 0 is nonempty, we
move a black ball from box 0 to box 1 with probability $(q - 1)^{-1}$.
The round halts when all the white balls are in box $F$.
\item[Round 2.] Note that, at the end of round 1, all the boxes are
empty but boxes 0 and 1 that contain only black
balls, and box $F$ that contains only white balls.
We paint in white all the balls in box 1 after which the game evolves
as described in round 1.
\end{longlist}
Any other round is defined starting from the final configuration of the
previous round in the same way as round 2 is defined starting from the
final configuration of round 1, and the game halts when all the balls
are either in box 0 or box $F$.
We refer the reader to the right-hand side of Figure \ref{figdomains}
for a schematic illustration of a single step evolution.
Note that, letting $\mathfrak A_j (t)$ denote the number of balls in
box $j$ at step $t$ for this game, it follows from Lemma \ref
{lemmoving} that
\[
\lim_{t \to\infty} E (\mathfrak A_0 (t)) \leq\lim_{t \to
\infty} E (\mathfrak B_0 (t)),
\]
whenever
\[
\mathfrak A_j (0) = \mathfrak B_j (0) \qquad\mbox
{for } j = 0, 1, \ldots, F;
\]
therefore it suffices to bound from below the limit on the left-hand side.
Let~$T_k$ denote the time at which round $k$ halts.
Since $F - j$ steps are required to move a white ball from box $j$ to
box $F$, and all the white balls are either in box 1 or box $F$ at the
beginning of round $k \geq2$,
\[
T_1 = \sum_{j = 1}^F (F - j) \mathfrak A_j (0) \quad\mbox
{and}\quad T_{k + 1} = T_k + (F - 1) \mathfrak A_1 (T_k).
\]
The expression of time $T_1$ together with Lemma \ref{leminitial}
implies that
\begin{eqnarray*}
E (T_1) & = &
\sum_{j = 1}^F (F - j) E (\mathfrak A_j (0)) = \sum_{j = 1}^F (F - j)
N p_j \\
&=&
\sum_{j = 0}^F (F - j) N p_j - N F p_0
\\
& = &
N F (1 - p_0) - N \sum_{j = 0}^F j p_j \\
&=&
N F \biggl(1 - \biggl(1 - \frac{1}{q} \biggr)^F -
\frac{1}{q} \biggr).
\end{eqnarray*}
In other respects, in view of the expression of time $T_{k + 1}$, and
since at each step a black ball is moved from box\vadjust{\goodbreak} 0 to box 1 with
probability $(q - 1)^{-1}$, we also have
\[
E (\mathfrak A_1 (T_{k + 1})) = \biggl(\frac{F - 1}{q - 1} \biggr) E (\mathfrak A_1
(T_k)) = \biggl(\frac{F - 1}{q - 1} \biggr)^k E (\mathfrak A_1 (T_1)).
\]
%
Combining the previous two equations, we obtain
\begin{eqnarray*}
E (\mathfrak A_1 (T_{k + 1})) &=& \biggl(\frac{F - 1}{q - 1} \biggr)^k \frac{1}{q -
1} E (T_1) \\
&=& N \biggl(\frac{F}{q - 1} \biggr) \biggl(\frac{F - 1}{q - 1}
\biggr)^k \biggl(1 - \biggl(1 - \frac{1}{q} \biggr)^F - \frac{1}{q}
\biggr).
\end{eqnarray*}
Finally, using again Lemma \ref{leminitial} and some basic algebra, we
deduce that
\begin{eqnarray*}
&&N^{-1} \lim_{t \to\infty} E (W_0 (t)) \\
&&\qquad \geq
N^{-1} \lim_{t \to\infty} E (\mathfrak A_0 (t)) \geq N^{-1} \Biggl(E
(\mathfrak A_0 (0)) - \sum_{k =
1}^{\infty} E (\mathfrak A_1 (T_k)) \Biggr) \\
&&\qquad\geq \biggl(1 - \frac{1}{q} \biggr)^F - \sum_{k = 0}^{\infty}
\biggl(\frac{F}{q - 1} \biggr) \biggl(\frac {F - 1}{q - 1} \biggr)^k
\biggl(1 - \biggl(1 - \frac{1}{q} \biggr)^F - \frac{1}{q} \biggr) \\
&&\qquad
\geq \biggl(1 - \frac{1}{q} \biggr)^F + \frac{F}{q - F}
\biggl(\biggl(1 - \frac{1}{q} \biggr)^F - \biggl(1 - \frac{1}{q}
\biggr) \biggr).
\end{eqnarray*}
This completes the proof.
\end{pf}

Theorem \ref{thcoexistence} directly follows from the combination of
Lemmas \ref{lemdomains} and \ref{lemfinal}.


%

%
\printaddresses


\begin{thebibliography}{14}

\bibitem{arratia1983}
%
\begin{barticle}[mr]
\bauthor{\bsnm{Arratia},~\bfnm{Richard}\binits{R.}}
(\byear{1983}).
\btitle{Site recurrence for annihilating random walks on $\mathbb Z^d$}.
\bjournal{Ann. Probab.}
\bvolume{11}
\bpages{706--713}.
\bid{issn={0091-1798}, mr={0704557}}
\bptok{imsref}%
\end{barticle}
%
\endbibitem

\bibitem{axelrod1997}
%
\begin{barticle}[auto:STB|2011/08/02|11:14:52]
\bauthor{\bsnm{Axelrod},~\bfnm{R.}\binits{R.}}
(\byear{1997}).
\btitle{The dissemination of culture: A model with local convergence
and global
polarization}.
\bjournal{J. Conflict Resolut.}
\bvolume{41}
\bpages{203--226}.
\bptok{imsref}%
\end{barticle}
%
\endbibitem

\bibitem{castellanofortunatoloreto2009}
%
\begin{barticle}[auto:STB|2011/08/02|11:14:52]
\bauthor{\bsnm{Castellano},~\bfnm{C.}\binits{C.}},
\bauthor{\bsnm{Fortunato},~\bfnm{S.}\binits{S.}} \AND
\bauthor{\bsnm{Loreto},~\bfnm{V.}\binits{V.}}
(\byear{2009}).
\btitle{Statistical physics of social dynamics}.
\bjournal{Rev. Modern Phys.}
\bvolume{81}
\bpages{591--646}.
\bptok{imsref}%
\end{barticle}
%
\endbibitem

\bibitem{cliffordsudbury1973}
%
\begin{barticle}[mr]
\bauthor{\bsnm{Clifford},~\bfnm{Peter}\binits{P.}} \AND
\bauthor{\bsnm{Sudbury},~\bfnm{Aidan}\binits{A.}}
(\byear{1973}).
\btitle{A model for spatial conflict}.
\bjournal{Biometrika}
\bvolume{60}
\bpages{581--588}.
\bid{issn={0006-3444}, mr={0343950}}
\bptok{imsref}%
\end{barticle}
%
\endbibitem

\bibitem{coxgriffeath1983}
%
\begin{barticle}[mr]
\bauthor{\bsnm{Cox},~\bfnm{J.~Theodore}\binits{J.~T.}} \AND
\bauthor{\bsnm{Griffeath},~\bfnm{David}\binits{D.}}
(\byear{1983}).
\btitle{Occupation time limit theorems for the voter model}.
\bjournal{Ann. Probab.}
\bvolume{11}
\bpages{876--893}.
\bid{issn={0091-1798}, mr={0714952}}
\bptok{imsref}%
\end{barticle}
%
\endbibitem

\bibitem{erdosney1974}
%
\begin{barticle}[mr]
\bauthor{\bsnm{Erd{\H{o}}s},~\bfnm{P.}\binits{P.}} \AND
\bauthor{\bsnm{Ney},~\bfnm{P.}\binits{P.}}
(\byear{1974}).
\btitle{Some problems on random intervals and annihilating particles}.
\bjournal{Ann. Probab.}
\bvolume{2}
\bpages{828--839}.
\bid{mr={0373068}}
\bptok{imsref}%
\end{barticle}
%
\endbibitem

\bibitem{harris1972}
%
\begin{barticle}[mr]
\bauthor{\bsnm{Harris},~\bfnm{T.~E.}\binits{T.~E.}}
(\byear{1972}).
\btitle{Nearest-neighbor {M}arkov interaction processes on multidimensional
lattices}.
\bjournal{Adv. Math.}
\bvolume{9}
\bpages{66--89}.
\bid{issn={0001-8708}, mr={0307392}}
\bptok{imsref}%
\end{barticle}
%
\endbibitem

\bibitem{holleyliggett1975}
%
\begin{barticle}[mr]
\bauthor{\bsnm{Holley},~\bfnm{Richard~A.}\binits{R.~A.}} \AND
\bauthor{\bsnm{Liggett},~\bfnm{Thomas~M.}\binits{T.~M.}}
(\byear{1975}).
\btitle{Ergodic theorems for weakly interacting infinite systems and
the voter
model}.
\bjournal{Ann. Probab.}
\bvolume{3}
\bpages{643--663}.
\bid{mr={0402985}}
\bptok{imsref}%
\end{barticle}
%
\endbibitem

\bibitem{lanchier2010}
%
\begin{barticle}[mr]
\bauthor{\bsnm{Lanchier},~\bfnm{N.}\binits{N.}}
(\byear{2010}).
\btitle{Opinion dynamics with confidence threshold: An alternative to the
{A}xelrod model}.
\bjournal{ALEA Lat. Am. J. Probab. Math. Stat.}
\bvolume{7}
\bpages{1--18}.
\bid{issn={1980-0436}, mr={2644039}}
\bptok{imsref}%
\end{barticle}
%
\endbibitem

\bibitem{schwartz1978}
%
\begin{barticle}[mr]
\bauthor{\bsnm{Schwartz},~\bfnm{Diane}\binits{D.}}
(\byear{1978}).
\btitle{On hitting probabilities for an annihilating particle model}.
\bjournal{Ann. Probab.}
\bvolume{6}
\bpages{398--403}.
\bid{mr={0494573}}
\bptok{imsref}%
\end{barticle}
%
\endbibitem

\bibitem{vazquezkrapivskyredner2003}
%
\begin{barticle}[mr]
\bauthor{\bsnm{Vazquez},~\bfnm{F.}\binits{F.}},
\bauthor{\bsnm{Krapivsky},~\bfnm{P.~L.}\binits{P.~L.}} \AND
\bauthor{\bsnm{Redner},~\bfnm{S.}\binits{S.}}
(\byear{2003}).
\btitle{Constrained opinion dynamics: Freezing and slow evolution}.
\bjournal{J. Phys. A}
\bvolume{36}
\bpages{L61--L68}.
\bid{doi={10.1088/0305-4470/36/3/103}, issn={0305-4470}, mr={1959419}}
\bptok{imsref}%
\end{barticle}
%
\endbibitem

\bibitem{vazquezredner2004}
%
\begin{barticle}[mr]
\bauthor{\bsnm{Vazquez},~\bfnm{F.}\binits{F.}} \AND
\bauthor{\bsnm{Redner},~\bfnm{S.}\binits{S.}}
(\byear{2004}).
\btitle{Ultimate fate of constrained voters}.
\bjournal{J. Phys. A}
\bvolume{37}
\bpages{8479--8494}.
\bid{doi={10.1088/0305-4470/37/35/006}, issn={0305-4470}, mr={2091256}}
\bptok{imsref}%
\end{barticle}
%
\endbibitem

\bibitem{vilonevespignanicastellano2002}
%
\begin{barticle}[auto:STB|2011/08/02|11:14:52]
\bauthor{\bsnm{Vilone},~\bfnm{D.}\binits{D.}},
\bauthor{\bsnm{Vespignani},~\bfnm{A.}\binits{A.}} \AND
\bauthor{\bsnm{Castellano},~\bfnm{C.}\binits{C.}}
(\byear{2002}).
\btitle{Ordering phase transition in the one-dimensional Axelrod model}.
\bjournal{Eur. Phys. J. B}
\bvolume{30}
\bpages{399--406}.
\bptok{imsref}%
\end{barticle}
%
\endbibitem

\end{thebibliography}
\end{document}